\documentclass[11pt,a4paper]{amsart}

\usepackage{amssymb,mathtools}
\usepackage[hidelinks]{hyperref}
\usepackage{microtype}
\hypersetup{
  pdftitle={Extreme Amenability and Failure of the Levy Property for the Trace-Class Special Unitary Group},
  pdfauthor={Sergio Luigi Cacciatori, Takashi Shioya, Pietro Ursino, Shigeaki Yokota},
  pdfsubject={Extreme amenability, non-Levy property, and trace-class unitary groups},
  pdfkeywords={extreme amenability, Levy groups, Schatten trace norm, special unitary group}
}

\newtheorem{thm}{Theorem}
\newtheorem{prop}[thm]{Proposition}
\newtheorem{lem}[thm]{Lemma}

\theoremstyle{remark}
\newtheorem{rem}[thm]{Remark}

\newcommand{\U}{\mathrm{U}}
\newcommand{\SU}{\mathrm{SU}}
\newcommand{\Tr}{\operatorname{Tr}}
\newcommand{\E}{\mathbb{E}}
\newcommand{\T}{\mathbb{T}}

\title[The trace-class special unitary group]{Extreme amenability and failure of the L\'evy property for the trace-class special unitary group}

\author[S. L. Cacciatori]{Sergio Luigi Cacciatori}
\address{Department of Science and High Technology, Universit\`a dell'Insubria,
Via Valleggio 11, IT-22100 Como, Italy; and INFN, Via Celoria 16,
IT-20133 Milano, Italy}
\email{sergio.cacciatori@uninsubria.it}

\author[T. Shioya]{Takashi Shioya}
\address{Mathematical Institute, Tohoku University, Sendai, Japan}
\email{shioya@math.tohoku.ac.jp}

\author[P. Ursino]{Pietro Ursino}
\address{Department of Mathematics and Informatics, Universit\`a degli Studi
 di Catania, Viale Andrea Doria 6, 95125 Catania, Italy; and National
 Group for Algebraic and Geometric Structures and their Applications
 (GNSAGA), INdAM, Italy}
\email{pietro.ursino@unict.it}

\author[S. Yokota]{Shigeaki Yokota}
\address{Department of Mathematics, Graduate School of Science,
The University of Osaka, Japan}
\email{shigeaki.yokota.t4@dc.tohoku.ac.jp}

\thanks{Takashi Shioya was supported by JSPS KAKENHI Grant Number 24K06729.
Pietro Ursino acknowledges financial support from the Universit\`a degli
Studi di Catania, PIAno di inCEntivi per la RIcerca di Ateneo 2024/2026 -
Linea di intervento 1.}



\date{}
\begin{document}
\begin{abstract}
We prove that the trace-class special unitary group, equipped with
the trace-norm topology, is extremely amenable but is not a L\'evy
group, answering questions posed by Pestov. The failure of the L\'evy
property follows from a uniform non-concentration estimate for every
nontrivial compact subgroup. Extreme amenability is established through
a transverse spectral-gap estimate and trace-norm duality, yielding
asymptotic invariance of Lipschitz observables under each fixed
finite-dimensional special unitary subgroup. A limiting argument then
gives fixed points for continuous actions on arbitrary compact Hausdorff
spaces. We also formulate an abstract criterion for extreme amenability
based on asymptotic quadratic invariance along a dense chain of compact
subgroups. Neither argument uses concentration results for other Schatten
exponents.
\end{abstract}
\maketitle

\section{Introduction}\label{sec:intro}

Extreme amenability is a strong fixed-point property in topological
dynamics. While amenability guarantees the existence of an invariant
probability measure on every nonempty compact space on which a group
acts continuously, extreme amenability requires a common fixed point.
Equivalently, the universal minimal flow of an extremely amenable group
is a singleton. The concept originated in the work of Mitchell and
Granirer in the 1960s; see \cite{Pestov}. Although no nontrivial locally
compact group is extremely amenable, Herer and Christensen
\cite{HererChristensen} constructed the first nontrivial example in 1975.
A major development was the work of Gromov and Milman
\cite{GromovMilman}, who established extreme amenability of the unitary
group of a separable infinite-dimensional Hilbert space with the strong
operator topology, using concentration of measure on high-dimensional
compact groups.

The concentration-of-measure approach of Gromov and Milman led to the
theory of \emph{L\'evy groups}. Their theorem asserts that every
L\'evy group is extremely amenable \cite{GromovMilman,Pestov}.
The converse, however, is false:
an extremely amenable group need not be a L\'evy group. In fact,
Farah and Solecki \cite{FarahSolecki} constructed an extremely amenable
group with a dense increasing union of compact subgroups which is not
a L\'evy group, thereby answering negatively an old question of
Furstenberg revisited by Pestov
\cite[\S5.2, Question~5, p.~214]{Pestov2005}.

Our work provides a new natural
example of this phenomenon in the nonabelian Banach--Lie group
$\SU_1(H)$,
together with a quantitatively controlled realization. Indeed, we found that concentration is impossible since the concentration function verifies $\alpha_K(1/4)\ge\frac{3}{128}$
for every nontrivial compact subgroup $K\subset\U_1(H)$.
Thus extreme amenability does not require concentration
of Haar measures along compact subgroup approximations. Another
important development was the connection with combinatorics established
by Kechris, Pestov, and Todor\v{c}evi\'c \cite{KPT}, who related extreme
amenability of automorphism groups of Fra\"iss\'e structures to structural
Ramsey theory.

We now consider the trace-class special unitary group.
Let $H$ be a separable infinite-dimensional complex Hilbert space.
Write $S_1(H)$ for the trace-class operators, with the
\emph{unnormalized} trace norm $\|A\|_{S_1}:=\Tr|A|$.
Define the trace-class unitary group and its determinant-one subgroup by
\[
\begin{aligned}
  \U_1(H)&:=\{U\in\U(H)\mid U-I\in S_1(H)\},\\
  \SU_1(H)&:=\ker\bigl(\det_F:\U_1(H)\longrightarrow\T\bigr),
\end{aligned}
\]
where $\det_F$ is the Fredholm determinant and
$\T:=\{z\in\mathbb C\mid |z|=1\}$.
Both groups carry the bi-invariant metric
\[
  d_1(U,V):=\|U-V\|_{S_1}.
\]
The distinction between $\U_1(H)$ and $\SU_1(H)$ is
essential. The Fredholm determinant is a continuous
surjective homomorphism from $\U_1(H)$ onto the circle
group. Consequently, $\U_1(H)$ is not extremely amenable.
This obstruction does not apply to $\SU_1(H)$, on which
the determinant is identically one.

The standard copies of $\SU(n)$, acting as the identity outside the
first $n$ coordinate vectors, have dense union in $\SU_1(H)$
(Lemma~\ref{lem:SU1-density}). This is the group denoted by
$\SU(\infty)_1$ in Pestov's terminology \cite[Remark~4.1.12]{Pestov}.

In his 2005 lecture notes
\cite[\S5.2, Question~7, pp.~214--215]{Pestov2005}, Pestov asked
whether $\SU_1(H)$ is a L\'evy
group and whether it is extremely amenable. The L\'evy-group
question also appears in his 2006 monograph
\cite[Remark~4.1.15]{Pestov}, alongside questions for other
Schatten-class unitary groups. We establish the following
contrasting conclusions.

\begin{thm}\label{thm:SU1-extreme}
The group $\SU_1(H)$ is extremely amenable with respect to the
trace-norm topology.
\end{thm}

\begin{thm}\label{thm:no-subgroup-levy}
No nontrivial subgroup of $\U_1(H)$, with the inherited trace-norm
topology, is a L\'evy group. In particular, $\SU_1(H)$ is not a
L\'evy group.
\end{thm}

Our results exhibit the distinction between extreme amenability and
the L\'evy property in a natural infinite-dimensional Banach--Lie
group. Although $\SU_1(H)$ contains a dense increasing union of compact
special unitary groups, no family of nontrivial compact subgroups can
form a L\'evy family for the trace metric. Nevertheless, $\SU_1(H)$ is
extremely amenable. In contrast to the classical Gromov--Milman
approach, the fixed-point property is obtained from asymptotic
invariance under each fixed finite-dimensional subgroup, rather than
concentration of the Haar measures on the entire approximating groups.
This provides a different mechanism for proving extreme amenability
and resolves the specific questions raised by Pestov
\cite{Pestov2005,Pestov}.

\paragraph{Relation to Schneider's concentration of invariant means.}
Schneider \cite{Schneider2023} proved a concentration inequality for
convolution products of invariant means adapted to chains of amenable
topological subgroups (Theorem~4.8), and obtained a sufficient criterion
for extreme amenability based on vanishing amenable length
(Theorem~4.12). His construction allows noncompact amenable subgroups,
whereas our argument uses the Haar probabilities of the compact chain
$(\SU(n))_{n\ge2}$. More importantly, our proof does not require these
Haar measures to concentrate: it relies instead on asymptotic
$L^2$-invariance of Lipschitz functions under each fixed $\SU(k)$,
derived from a relative spectral gap and trace-norm duality.
There is also a concrete distinction between the two sufficient
criteria. Schneider showed that every nontrivial topological group
with the no-small-subgroups property has amenable length at least one
\cite[Remark~4.11(4)]{Schneider2023}. Since $\SU_1(H)$ is a
Banach--Lie group and has the no-small-subgroups property, the
vanishing-amenable-length condition of Schneider's Theorem~4.12
cannot be applied directly in the present setting. Our relative
spectral-energy criterion, Proposition~\ref{prop:abstract-spectral},
provides a different route to extreme amenability.

Theorem~\ref{thm:no-subgroup-levy} is proved using a uniform
quantitative obstruction
of the concentration function as stated before.
A compact subgroup
necessarily acts trivially off a finite-dimensional subspace, and
randomizing trace observables on this subspace gives the bound.
This rules out all compact-subgroup approximations, not just the
standard chain $\SU(n)$.

The proof of Theorem~\ref{thm:SU1-extreme} is independent of that
obstruction. For $2\le k<n$, if $f:\SU(n)\to\mathbb R$ is
$L$-Lipschitz for the trace metric, we show that
\begin{equation}\label{eq:intro-translation}
  \sup_{h\in\SU(k)}\int_{\SU(n)}
  |f(hu)-f(u)|^2\,d\nu_n(u)
  \le\frac{8kL^2}{n-k},
\end{equation}
where $\nu_n$ is the Haar probability measure. The argument combines a
transverse relative spectral gap, obtained from highest weights, with
a dimension-free energy estimate from trace-norm duality.
For a continuous action on a compact Hausdorff space,
the resulting asymptotic $L^2$-invariance implies that
every weak-$^*$ cluster point of the orbit measures is
supported on the common fixed-point set. This yields
extreme amenability. We first formulate the fixed-point mechanism as an abstract
asymptotic-invariance criterion (Theorem~\ref{thm:abstract-extreme}).
From it we deduce a relative spectral-energy criterion
(Proposition~\ref{prop:abstract-spectral}), which is then applied
directly to prove Theorem~\ref{thm:SU1-extreme}.
The resulting principle applies to topological groups beyond the
present trace-class setting, provided the quantitative hypotheses
can be verified.

For comparison, let $1\le p<\infty$ and consider the
Schatten-class unitary group
\[
  \U_p(H):=\{U\in\U(H)\mid U-I\in S_p(H)\},
\]
where $S_p(H)$ is the Schatten $p$-class. We equip
$\U_p(H)$ with the unnormalized Schatten $p$-norm topology.
Let $\SU_p(H)$ denote the closure of $\bigcup_{n\ge2}\SU(n)$
in $\U_p(H)$. As observed by Pestov
\cite[Remark~4.1.12]{Pestov}, $\SU_p(H)=\U_p(H)$ for $p>1$,
whereas $\SU_1(H)$ is a proper subgroup of $\U_1(H)$.
The cases $p=1$ and $p=2$ were already settled, as recorded by
Pestov \cite[Remark~4.1.12 and Theorem~4.1.13]{Pestov}:
$\U_1(H)$ is not a L\'evy group, whereas $\U_2(H)$ is a
L\'evy group.

\begin{prop}\label{prop:Up-levy}
For $1\le p<\infty$, the group $\U_p(H)$ is a L\'evy group
if and only if $p>1$.
\end{prop}

We prove Proposition~\ref{prop:Up-levy} in
Section~\ref{sec:schatten-levy} using classical concentration
results. Its proof is independent of those of our main theorems.

\section{Preliminaries and finite-dimensional approximation}
\label{sec:prelim}

A \emph{metric probability space} is a complete separable metric space equipped with a Borel probability measure.
For a metric probability space $X$, we denote its metric and measure by $d_X$ and $\mu_X$, respectively.

Let $X$ be a metric probability space.
For a Borel set $B\subset X$, set
$B_r:=\{x\in X\mid d_X(x,B)<r\}$, where $d_X(x,B):= \inf_{b \in B} d_X(x,b)$.
We define the \emph{concentration function} of $X$ by
\begin{equation}\label{eq:concentration-definition}
  \alpha_X(r):=\sup_{\mu_X(B)\ge1/2}\bigl(1-\mu_X(B_r)\bigr),\qquad r>0.
\end{equation}
A directed family of metric probability spaces is a \emph{L\'evy family} if
its concentration functions tend to zero at every positive radius.
A topological group with a compatible bi-invariant metric is a
\emph{L\'evy group} if it has a directed increasing family of
compact subgroups with dense union whose normalized Haar probability
spaces form a L\'evy family for the restricted metric \cite{Pestov}.
A group is \emph{extremely amenable} if every continuous action on a
nonempty compact Hausdorff space has a fixed point.

A \emph{median} of a real-valued Borel function $f$ on $X$ is a real number $m$
with $\mu_X\{f\le m\}\ge1/2$ and $\mu_X\{f\ge m\}\ge1/2$.
For every real-valued $1$-Lipschitz function $f$ on $X$ we have
\begin{equation}\label{eq:median-tail}
  \mu_X\{|f-m|\ge r\}\le2\alpha_X(r),\qquad r>0.
\end{equation}
Indeed, apply the definition to the $r$-neighborhoods of the two
sets $\{f\le m\}$ and $\{f\ge m\}$.

Let $H$ be a separable infinite-dimensional complex Hilbert space and
let $\mathcal{B}(H)$ denote the Banach algebra of all bounded linear
operators on $H$, equipped with the operator norm
\[
  \|T\|_{\mathrm{op}}:=\sup_{\|x\|_H\le 1}\|Tx\|_H.
\]
The trace-class operators embed continuously into $\mathcal{B}(H)$,
since
\[
  \|T\|_{\mathrm{op}}\le \|T\|_{S_1}
  \qquad (T\in S_1(H)).
\]

Fix an orthonormal basis of $H$. The standard copies of $\U(n)$ and
$\SU(n)$ act as the identity outside its first $n$ vectors. For $n\ge2$,
write
\[
  d_{1,n}(u,v):=\|u-v\|_{S_1},\qquad u,v\in\SU(n),
\]
for the restriction of the unnormalized trace metric. For the rest of
the paper, $\nu_n$ denotes normalized Haar measure on $\SU(n)$.

\begin{lem}\label{lem:SU1-density}
The union $\bigcup_{n\ge2}\SU(n)$ is dense in $\SU_1(H)$ for $d_1$.
In particular, $\SU_1(H)$ is its trace-norm completion.
\end{lem}

\begin{proof}
The density assertion follows from
\cite[Remark~4.1.12]{Pestov}.
The completeness follows from the completeness of
$\U_1(H)$ in the trace norm and the fact that
$\SU_1(H)=\ker\det_F$ is closed.
\end{proof}

\section{Compact subgroups and non-concentration}\label{sec:nonlevy}

Throughout this section compactness is taken in the trace-norm topology.

\subsection{Finite-dimensional support of compact subgroups}

For a compact subgroup $K\subset\U_1(H)$, let $\mu_K$ denote its
normalized Haar measure and set
\[
  H^K:=\{\xi\in H\mid U\xi=\xi\text{ for every }U\in K\},
  \qquad E_K:=(H^K)^\perp.
\]

\begin{lem}\label{lem:compact-support}
For every compact subgroup $K\subset\U_1(H)$, the subspace $E_K$ is
finite-dimensional and $K$-invariant. Every $U\in K$ acts as the
identity on $E_K^\perp$, and
\begin{equation}\label{eq:mean-representation-zero}
  \int_K U|_{E_K}\,d\mu_K(U)=0.
\end{equation}
If $K$ is nontrivial, then $\dim E_K\ge1$.
\end{lem}

\begin{proof}
The continuous map $U\mapsto I-U$ from $K$ to the trace-class Banach
space $S_1(H)$ is Bochner integrable. Put
\[
  Q:=\int_K(I-U)\,d\mu_K(U)\in S_1(H),\qquad P:=I-Q.
\]
The inclusion $S_1(H)\hookrightarrow\mathcal{B}(H)$ is continuous,
so for each $\xi\in H$,
\[
  P\xi=\int_K U\xi\,d\mu_K(U).
\]
Haar invariance implies that $P\xi\in H^K$, and $P$ is the identity
on $H^K$. Moreover, invariance under inversion and $U^*=U^{-1}$ show
that $P$ is self-adjoint. Thus $P$ is the orthogonal projection onto
$H^K$, and $Q=I-P$ is the orthogonal projection onto $E_K$.
A trace-class orthogonal projection has finite rank, since each nonzero
singular value equals one. Hence $E_K$ is finite-dimensional.
The remaining assertions follow from the definition of $H^K$ and
$P|_{E_K}=0$. In particular, $E_K=\{0\}$ would force $K=\{I\}$.
\end{proof}

\subsection{Low-order Haar moments}

We record the elementary moments used in the uniform obstruction.
They are special cases of the unitary trace moment calculations in
\cite{DiaconisShahshahani}; the short proof below also covers the
small dimensions explicitly.

\begin{lem}\label{lem:Haar-moments}
For a Haar-distributed matrix $W\in\U(N)$, $N\ge1$,
\begin{align}
  \E|\Tr W|^2&=1,\label{eq:complex-second-moment}\\
  \E|\Tr W|^4&=
\begin{cases}
    1,&N=1,\\
    2,&N\ge2.
\end{cases}\label{eq:complex-fourth-moment}
\end{align}
Consequently,
\begin{equation}\label{eq:real-Haar-moments}
  \E(\operatorname{Re}\Tr W)^2=\frac{1}{2},
  \qquad
  \E(\operatorname{Re}\Tr W)^4\le\frac{3}{4}.
\end{equation}
\end{lem}

\begin{proof}
The defining representation of $\U(N)$ is irreducible, so the Schur orthogonality relation for its character $\chi(W):=\Tr W$ gives
\eqref{eq:complex-second-moment}. If $N\ge2$, its tensor square is the
sum of two nonisomorphic irreducible representations,
\[
  \mathbb{C}^N\otimes\mathbb{C}^N
  =\operatorname{Sym}^2(\mathbb{C}^N)
   \oplus\bigwedge\nolimits^2(\mathbb{C}^N).
\]
Since the character of this tensor square is $\chi^2$, character
orthogonality gives $\E|\chi(W)|^4=2$. For $N=1$ this moment is one,
since $|\Tr W|=1$.

The distribution of $\Tr W$ is invariant under multiplication by any
scalar of modulus one. Averaging over this scalar and using
\[
  \frac{1}{2\pi}\int_0^{2\pi}\cos^2\theta\,d\theta=\frac{1}{2},
  \qquad
  \frac{1}{2\pi}\int_0^{2\pi}\cos^4\theta\,d\theta=\frac{3}{8}
\]
gives
\[
  \E(\operatorname{Re}\Tr W)^2=\frac{1}{2}\E|\Tr W|^2,
  \qquad
  \E(\operatorname{Re}\Tr W)^4=\frac{3}{8}\E|\Tr W|^4.
\]
This proves \eqref{eq:real-Haar-moments}.
\end{proof}

\subsection{A uniform non-concentration bound}

For a compact subgroup $K\subset\U_1(H)$, write $\alpha_K$ for the
concentration function of $K$ with the restricted trace metric $d_1$
and normalized Haar measure $\mu_K$.

\begin{thm}\label{thm:compact-uniform}
Every nontrivial compact subgroup $K\subset\U_1(H)$ satisfies
\begin{equation}\label{eq:uniform-nonconcentration}
  \alpha_K(1/4)\ge\frac{3}{128}.
\end{equation}
\end{thm}

\begin{proof}
By Lemma~\ref{lem:compact-support}, $E:=E_K$ has finite dimension
$N\ge1$. For $U\in K$, write $U_E:=U|_E$. Since all elements of $K$
are the identity on $E^\perp$,
\begin{equation}\label{eq:restricted-trace-metric}
  d_1(U,V)=\|U_E-V_E\|_{S_1(E)},\qquad U,V\in K.
\end{equation}
Here $S_1(E)$ denotes the trace class on $E$ with its unnormalized norm,
and $\Tr_E$ denotes the trace on $E$. Let $\eta_E$ be normalized Haar
measure on $\U(E)$.

For $A\in\U(E)$ define the real-valued function
\[
  f_A(U):=\operatorname{Re}\Tr_E(AU_E),\qquad U\in K.
\]
Trace duality and \eqref{eq:restricted-trace-metric} give
\begin{align*}
  |f_A(U)-f_A(V)|
  &\le |\Tr_E(A(U_E-V_E))|\\
  &\le\|A\|_{\mathrm{op}}\|U_E-V_E\|_{S_1(E)}
   =d_1(U,V).
\end{align*}
Thus $f_A$ is $1$-Lipschitz. By
\eqref{eq:mean-representation-zero}, it is centered for every $A$:
\begin{equation}\label{eq:randomized-trace-centered}
  \int_K f_A(U)\,d\mu_K(U)=0.
\end{equation}
The auxiliary unitary $A$ need not belong to the image of $K$ in $\U(E)$.

Set
\[
  v(A):=\int_K f_A(U)^2\,d\mu_K(U),\qquad
  q(A):=\int_K f_A(U)^4\,d\mu_K(U).
\]
For fixed $U\in K$, the matrix $AU_E$ is Haar distributed on $\U(E)$
when $A$ has distribution $\eta_E$. Fubini's theorem and
Lemma~\ref{lem:Haar-moments} therefore yield
\begin{equation}\label{eq:averaged-moments}
  \E_Av(A)=\frac{1}{2},\qquad \E_Aq(A)\le\frac{3}{4},
\end{equation}
where $\E_A$ denotes integration with respect to $\eta_E$.

Now choose $A,U,V$ independently with distributions
$\eta_E,\mu_K,\mu_K$, respectively, and put
\[
  D:=f_A(U)-f_A(V).
\]
Conditional on $A$, equation~\eqref{eq:randomized-trace-centered} and
independence give
\[
  \E_{U,V}D^2=2v(A),\qquad
  \E_{U,V}D^4=2q(A)+6v(A)^2.
\]
Since $v(A)^2\le q(A)$ by Cauchy--Schwarz, averaging in $A$ gives
\begin{equation}\label{eq:difference-moments}
  \E D^2=1,\qquad \E D^4\le8\E_Aq(A)\le6.
\end{equation}
A further application of Cauchy--Schwarz gives
\begin{align*}
  1=\E D^2
  &\le\frac{1}{4}+
       \E\bigl[D^2\mathbf{1}_{\{|D|>1/2\}}\bigr]\\
  &\le\frac{1}{4}+
       \sqrt{\E D^4\,\mathbb{P}(|D|>1/2)}\\
  &\le\frac{1}{4}+\sqrt{6\,\mathbb{P}(|D|>1/2)}.
\end{align*}
Consequently,
\begin{equation}\label{eq:difference-lower-tail}
  \mathbb{P}(|D|>1/2)\ge\frac{3}{32}.
\end{equation}

On the other hand, fix $A$ and choose a median $m_A$ of $f_A$.
By \eqref{eq:median-tail},
\[
  \mu_K\{|f_A-m_A|>r\}\le2\alpha_K(r).
\]
If $|f_A(U)-f_A(V)|>2r$, then at least one of
$|f_A(U)-m_A|$ and $|f_A(V)-m_A|$ exceeds $r$. Hence
\begin{equation}\label{eq:difference-upper-tail}
  \mathbb{P}_{U,V}\bigl(|f_A(U)-f_A(V)|>2r\bigr)
  \le4\alpha_K(r).
\end{equation}
This inequality holds for every fixed $A$, so it can be integrated with
respect to $\eta_E$.
Taking $r=1/4$ and combining
\eqref{eq:difference-lower-tail} with \eqref{eq:difference-upper-tail},
we obtain
\[
  \frac{3}{32}\le\mathbb{P}(|D|>1/2)\le4\alpha_K(1/4),
\]
which proves \eqref{eq:uniform-nonconcentration}.
\end{proof}

\begin{proof}[Proof of Theorem~\ref{thm:no-subgroup-levy}]
Let $G\subset\U_1(H)$ be a nontrivial subgroup with the inherited
trace-norm topology. The restriction of $d_1$ is a compatible
bi-invariant metric on $G$. Suppose that $(K_i)_{i\in I}$ is an
increasing family of compact subgroups with dense union witnessing
the L\'evy property of $G$. Compactness in $G$ implies compactness
in $\U_1(H)$, because the inclusion is continuous.
Since $G$ is nontrivial and $\{I\}$ is closed, some $K_{i_0}$ is
nontrivial. Then every $K_i$ with $i\ge i_0$ is nontrivial, so
Theorem~\ref{thm:compact-uniform} gives
\[
  \alpha_{K_i}(1/4)\ge\frac{3}{128}\qquad(i\ge i_0).
\]
Thus $G$ is not a L\'evy group.
In particular, $\SU_1(H)$ is nontrivial, since it contains the
standard copy of $\SU(2)$, and the assertion applies to it.
\end{proof}

\begin{rem}
By Theorem~\ref{thm:compact-uniform}, any family consisting of
nontrivial compact subgroups of $\U_1(H)$, with their restricted
trace metrics and normalized Haar measures, fails to be a L\'evy
family, whether or not it is increasing. In particular,
$\alpha_{\SU(n)}(1/4)\ge3/128$ for every $n\ge2$. This conclusion
does not require a central limit theorem for traces. The estimate
does not rule out extreme amenability, which is proved independently
in Section~\ref{sec:extreme}.
\end{rem}

\section{Extreme amenability via transverse spectral estimates}
\label{sec:extreme}

We now prove Theorem~\ref{thm:SU1-extreme}. The uniform
non-concentration bound in Theorem~\ref{thm:compact-uniform} does not
preclude the fixed-point property. The relevant substitute for ordinary
concentration is \emph{asymptotic invariance under each fixed
finite-dimensional subgroup}. We obtain it from a transverse
spectral gap, together with trace-norm duality.

Throughout this section, $\nu_n$ denotes normalized Haar measure on
$\SU(n)$. For $2\le k<n$, let $\SU(k)$ be the standard subgroup of $\SU(n)$ acting
on the first $k$ coordinates, and put
\[
  (P_kf)(u):=\int_{\SU(k)}f(h^{-1}u)\,d\nu_k(h),
  \qquad f\in L^2(\SU(n),\nu_n).
\]
Then $P_k$ is the orthogonal projection onto the subspace of functions
invariant under left $\SU(k)$-translations.

\subsection{Transverse spectral and energy estimates}

Consider the orthogonal decomposition
$\mathbb C^n=\mathbb C^k\oplus\mathbb C^{n-k}$.
The space of \emph{transverse directions} is defined by
\[
\mathfrak m_{k,n}:=
\left\{
\begin{pmatrix}
0&B\\
-B^*&0
\end{pmatrix}
\;\middle|\;
B\in M_{k,n-k}(\mathbb C)
\right\}
\subset\mathfrak{su}(n).
\]
Let $E_{ij}$ be the matrix units and, for $1\le i\le k<j\le n$, set
\[
  X_{ij}:=\frac{E_{ij}-E_{ji}}{\sqrt2},
  \qquad
  Y_{ij}:=\frac{\mathrm i(E_{ij}+E_{ji})}{\sqrt2}.
\]
The matrices $X_{ij}$ and $Y_{ij}$ form a Hilbert--Schmidt
orthonormal basis of $\mathfrak m_{k,n}$. For $X\in\mathfrak{su}(n)$, write
\[
  (D_Xf)(u):=\left.\frac{d}{dt}\right|_{t=0} f(e^{tX}u).
\]
Let $\lambda$ denote the left regular representation,
\[
  (\lambda(a)f)(u):=f(a^{-1}u),\qquad a,u\in\SU(n).
\]
Thus $d\lambda(X)=-D_X$ on smooth functions. This sign disappears
in the squares defining the associated nonnegative operator.
For a smooth real- or complex-valued function $f$ define
\begin{equation}\label{eq:cross-energy}
  \mathcal E_{k,n}(f):=
  \sum_{i\le k<j}\left(
     \|D_{X_{ij}}f\|_{L^2(\nu_n)}^2+
     \|D_{Y_{ij}}f\|_{L^2(\nu_n)}^2
  \right).
\end{equation}

\begin{lem}\label{lem:relative-gap}
For every $2\le k<n$ and every smooth
$f:\SU(n)\to\mathbb C$,
\begin{equation}\label{eq:relative-gap}
  (n-k)\|f-P_kf\|_{L^2(\nu_n)}^2
  \le\mathcal E_{k,n}(f).
\end{equation}
\end{lem}

\begin{proof}
We first prove the corresponding operator inequality in an arbitrary
finite-dimensional unitary representation
$\pi:\SU(n)\to\U(V)$. We then apply it to the irreducible summands
of the left regular representation $\lambda$ and use the
Peter--Weyl decomposition.

Extend $d\pi$ complex linearly to $\mathfrak{sl}_n(\mathbb C)$ and
write $T_{ij}:=d\pi(E_{ij})$ for $i\ne j$. Unitarity implies
$T_{ij}^*=T_{ji}$. The nonnegative operator associated with the
energy in \eqref{eq:cross-energy} on $V$ is
\begin{equation}\label{eq:cross-Casimir}
\begin{aligned}
  Q_{k,n}
  &:=-\sum_{i\le k<j}
    \bigl(d\pi(X_{ij})^2+d\pi(Y_{ij})^2\bigr)\\
  &=\sum_{i\le k<j}
    \bigl(T_{ij}T_{ji}+T_{ji}T_{ij}\bigr).
\end{aligned}
\end{equation}
The span of the $X_{ij}$ and $Y_{ij}$ is invariant under the
adjoint action of $\SU(k)$. Hence $Q_{k,n}$ commutes with $\pi(\SU(k))$.

Consider a $Q_{k,n}$-eigenspace with eigenvalue $q$ that contains a
nontrivial irreducible $\SU(k)$-subrepresentation. Since the
eigenspace is $\SU(k)$-invariant, we can choose a unit highest-weight
vector $v$ of such a subrepresentation within it. For the standard
positive roots of $\mathfrak{sl}_k(\mathbb C)$, its highest weight is
nonzero, so
\begin{equation}\label{eq:highest-root-value}
  d\pi(E_{11}-E_{kk})v=\ell v,
  \qquad \ell\in\mathbb Z,\quad\ell\ge1.
\end{equation}
To make the last inequality explicit, let $\Lambda$ be this highest
weight and put
\[
  m_r:=\Lambda(E_{rr}-E_{r+1,r+1}),\qquad 1\le r<k.
\]
The $m_r$ are nonnegative integers and are not all zero. Since
$E_{11}-E_{kk}=\sum_{r=1}^{k-1}(E_{rr}-E_{r+1,r+1})$,
we have $\ell=\sum_{r=1}^{k-1}m_r\ge1$.

For each $j>k$, the commutator identities
\[
  [T_{ij},T_{ji}]=d\pi(E_{ii}-E_{jj})
  \qquad(i=1,k)
\]
and the inequalities $a+b\ge|a-b|$ for $a,b\ge0$ give
\begin{align*}
 &\sum_{i\in\{1,k\}}
       \bigl(\|T_{ij}v\|^2+\|T_{ji}v\|^2\bigr)\\
 &\quad\ge
  \left|\left\langle d\pi(E_{11}-E_{jj})v,v\right\rangle\right|
  +\left|\left\langle d\pi(E_{kk}-E_{jj})v,v\right\rangle\right|\\
 &\quad\ge
  \left|\left\langle d\pi(E_{11}-E_{kk})v,v\right\rangle\right|
  =\ell\ge1.
\end{align*}
Summing over $j=k+1,\ldots,n$ and using
\eqref{eq:cross-Casimir}, we obtain
\[
  q=\langle Q_{k,n}v,v\rangle\ge(n-k)\ell\ge n-k.
\]
Thus every eigenspace of $Q_{k,n}$ having a nontrivial
$\SU(k)$-type has eigenvalue at least $n-k$. On the $\SU(k)$-fixed
subspace $Q_{k,n}$ is merely nonnegative. If $P_K$ denotes the
orthogonal projection onto this fixed subspace, it follows that
\[
  Q_{k,n}\ge(n-k)(I-P_K).
\]
To pass to $L^2(\SU(n),\nu_n)$, write
$f=\sum_{\pi\in\widehat{\SU(n)}}f_\pi$ for its Peter--Weyl
decomposition into $\lambda$-isotypic components. Here
$\widehat{\SU(n)}$ denotes the set of equivalence classes of
irreducible unitary representations. Each isotypic component is a
finite orthogonal sum of copies of the corresponding representation.
The component projections commute with left translations, with
$P_k$, and, on smooth functions, with every $D_X$. The operator
inequality above therefore applies on each component. For smooth
$f$, Parseval's identity yields
\begin{align*}
  \|f-P_kf\|_{L^2(\nu_n)}^2
    &=\sum_{\pi\in\widehat{\SU(n)}}
        \|f_\pi-P_kf_\pi\|_{L^2(\nu_n)}^2,\\
  \mathcal E_{k,n}(f)
    &=\sum_{\pi\in\widehat{\SU(n)}}
        \mathcal E_{k,n}(f_\pi).
\end{align*}
All these sums converge, since $f$ and the finitely many derivatives
occurring in \eqref{eq:cross-energy} belong to $L^2(\nu_n)$.
Summing the component inequalities proves \eqref{eq:relative-gap}.
\end{proof}

The key feature of the next estimate is that its right-hand side
depends on $k$, but not on $n$.

\begin{lem}\label{lem:cross-energy-bound}
If $f:\SU(n)\to\mathbb R$ is smooth and $L$-Lipschitz for the
unnormalized trace metric $d_{1,n}$, then
\begin{equation}\label{eq:energy-bound}
  \mathcal E_{k,n}(f)\le2kL^2.
\end{equation}
\end{lem}

\begin{proof}
Fix $u\in\SU(n)$ and consider the real linear functional on
traceless Hermitian matrices
\[
  A\longmapsto D_{\mathrm iA}f(u).
\]
The $d_{1,n}$-Lipschitz condition implies
\[
  |D_{\mathrm iA}f(u)|\le L\|A\|_{S_1}.
\]
By the real Hahn--Banach theorem and the duality between the trace
and operator norms, there is a Hermitian matrix $B$ (depending on $u$) with
$\|B\|_{\mathrm{op}}\le L$ and
\begin{equation}\label{eq:dual-derivative}
  D_{\mathrm iA}f(u)=\Tr(BA)
\end{equation}
for every traceless Hermitian $A$. No measurable choice of $B(u)$
is required, because the following bound is pointwise.

Applying \eqref{eq:dual-derivative} to $A=-\mathrm iX_{ij}$ and
$A=-\mathrm iY_{ij}$ gives
\[
  |D_{X_{ij}}f(u)|^2+|D_{Y_{ij}}f(u)|^2
  =2|B_{ij}|^2.
\]
Consequently,
\begin{align*}
 &\sum_{i\le k<j}
    \bigl(|D_{X_{ij}}f(u)|^2+|D_{Y_{ij}}f(u)|^2\bigr)\\
 &\qquad=2\sum_{i\le k<j}|B_{ij}|^2
 \le 2\sum_{i=1}^k\sum_{j=1}^n|B_{ij}|^2\\
 &\qquad=2\sum_{i=1}^k\langle B^2e_i,e_i\rangle
 \le 2kL^2.
\end{align*}
Integration over $u$ proves \eqref{eq:energy-bound}.
\end{proof}

\subsection{An abstract asymptotic-invariance criterion}
\label{subsec:abstract-criterion}

  We first isolate an abstract fixed-point principle for
asymptotic quadratic invariance. We then derive from it the
spectral-energy criterion used directly in the proof of
Theorem~\ref{thm:SU1-extreme}. Neither criterion requires global
concentration of Haar measure or a L\'evy family.

For a real-valued Lipschitz function $f$ on a metric subspace
$Y\subseteq G$, define its Lipschitz constant by
\[
  \operatorname{Lip}_d(f):=
  \sup_{\substack{x,y\in Y\\x\ne y}}
  \frac{|f(x)-f(y)|}{d(x,y)}.
\]

\begin{thm}[Asymptotic-invariance criterion]
\label{thm:abstract-extreme}
Let $G$ be a topological group equipped with a compatible
right-invariant metric $d$. Suppose that
\[
  K_2\subseteq K_3\subseteq\cdots\subseteq G
\]
is an increasing sequence of compact subgroups with dense union,
and denote normalized Haar measure on $K_n$ by $\mu_n$.
Assume that for each fixed $k\ge2$ there are nonnegative numbers
$\varepsilon_{k,n}$, $n>k$, such that
$\varepsilon_{k,n}\to0$ as $n\to\infty$ and, for every
real-valued Lipschitz function $f$ on $K_n$,
\begin{equation}\label{eq:abstract-L2}
  \sup_{h\in K_k}\int_{K_n}|f(hu)-f(u)|^2\,d\mu_n(u)
  \le \varepsilon_{k,n}\operatorname{Lip}_d(f)^2.
\end{equation}
Then $G$ is extremely amenable.
\end{thm}

\begin{proof}
 Let $G$ act continuously on a nonempty compact Hausdorff space $X$,
and fix $x_0\in X$. Given $\varphi\in C(X,\mathbb R)$, define
  $F_\varphi(g):=\varphi(gx_0)$. This function is bounded. Moreover,
 continuity of the action and compactness of $X$ imply
\[
  \sup_{x\in X}|\varphi(ax)-\varphi(x)|\longrightarrow0
  \qquad(a\to e).
\]
  Since $d$ is right-invariant,
$d(g',g)=d(g'g^{-1},e)$; hence
\[
  |F_\varphi(g')-F_\varphi(g)|\le
  \sup_{x\in X}|\varphi(g'g^{-1}x)-\varphi(x)|.
\]
Thus $F_\varphi$ is uniformly continuous for $d$.

  For any $\eta>0$, bounded uniform continuity permits a bounded
$L$-Lipschitz function $F_L:G\to\mathbb R$ with
$\|F_\varphi-F_L\|_\infty<\eta$; for instance, take the infimal
convolution
\[
  F_L(g):=\inf_{v\in G}\bigl(F_\varphi(v)+Ld(g,v)\bigr)
\]
  with $L$ sufficiently large. Indeed, if
$|F_\varphi(g)-F_\varphi(v)|<\eta/2$ whenever $d(g,v)<\delta$,
 choose $L\delta>\sup F_\varphi-\inf F_\varphi$.
The infimum is then attained, up to arbitrarily small errors,
within distance $\delta$ of $g$, and
$0\le F_\varphi-F_L\le\eta/2$.

Fix $k\ge2$ and $h\in K_k$. For every $n>k$, the restriction of
$F_L$ to $K_n$ is $L$-Lipschitz, so
\begin{align*}
 &\left(\int_{K_n}
    |\varphi(hux_0)-\varphi(ux_0)|^2\,d\mu_n(u)\right)^{1/2}\\
 &\qquad\le 2\eta+L\sqrt{\varepsilon_{k,n}}.
\end{align*}
  The constant $L$ is independent of $n$. Therefore, letting
 $n\to\infty$ and then $\eta\downarrow0$, we find
\begin{equation}\label{eq:abstract-orbit}
  \lim_{n\to\infty}\int_{K_n}
  |\varphi(hux_0)-\varphi(ux_0)|^2\,d\mu_n(u)=0
\end{equation}
for every $h\in\bigcup_{k\ge2}K_k$ and
 $\varphi\in C(X,\mathbb R)$.

Push $\mu_n$ forward by $u\mapsto ux_0$ to obtain
$\rho_n\in\mathcal P(X)$. Compactness of the space of regular Borel
probability measures on $X$ in the weak-$^*$ topology gives a
convergent subnet $\rho_{n_j}\to\rho$.
  As the original index sequence tends to infinity along this
subnet, equation~\eqref{eq:abstract-orbit} implies
\[
  \int_X|\varphi(hx)-\varphi(x)|^2\,d\rho(x)=0.
\]
For each fixed $h$ and $\varphi$, the integrand is continuous and
nonnegative, so it vanishes at every point of
 $\operatorname{supp}\rho$. Since continuous real functions separate
points of a compact Hausdorff space, every point of
$\operatorname{supp}\rho$ is fixed by every
$h\in\bigcup_{k\ge2}K_k$ and hence, by density and continuity,
by all of $G$. The support is nonempty, giving a fixed point.
\end{proof}

The following proposition is a direct consequence of
Theorem~\ref{thm:abstract-extreme}. It converts a relative
spectral gap and a Lipschitz-energy bound into the asymptotic
invariance required by that theorem.

\begin{prop}[Relative spectral-energy criterion]
\label{prop:abstract-spectral}
Under the structural assumptions on $G$, $d$ and $(K_n)$ in
Theorem~\ref{thm:abstract-extreme}, define
\[
  (P_{k,n}f)(u):=\int_{K_k} f(h^{-1}u)\,d\mu_k(h)
\]
on $L^2(K_n,\mu_n)$. Suppose that for $2\le k<n$ there are
positive numbers $\lambda_{k,n}$ and nonnegative numbers
$C_{k,n}$, together with nonnegative functionals
$\mathcal E_{k,n}$ on real-valued Lipschitz functions on $K_n$,
such that
\begin{align}
  \|f-P_{k,n}f\|_{L^2(\mu_n)}^2
    &\le \lambda_{k,n}^{-1}\mathcal E_{k,n}(f),
\label{eq:abstract-gap}\\
  \mathcal E_{k,n}(f)
    &\le C_{k,n}\operatorname{Lip}_d(f)^2.
\label{eq:abstract-energy}
\end{align}
If $C_{k,n}/\lambda_{k,n}\to0$ as $n\to\infty$ for each
fixed $k$, then $G$ is extremely amenable.
\end{prop}

\begin{proof}
Write $(L_hf)(u):=f(hu)$ for $h\in K_k$.
By Haar invariance, $L_h$ is unitary on $L^2(K_n,\mu_n)$;
moreover, $L_hP_{k,n}f=P_{k,n}f$.
Consequently,
\[
 \|L_hf-f\|_{L^2(\mu_n)}
 \le \|L_h(f-P_{k,n}f)\|_2+\|P_{k,n}f-f\|_2
 =2\|f-P_{k,n}f\|_2.
\]
Combining this inequality with \eqref{eq:abstract-gap} and
\eqref{eq:abstract-energy} yields
\[
  \sup_{h\in K_k}\int_{K_n}|f(hu)-f(u)|^2\,d\mu_n(u)
  \le \frac{4C_{k,n}}{\lambda_{k,n}}
        \operatorname{Lip}_d(f)^2.
\]
For each fixed $k$, the ratio $C_{k,n}/\lambda_{k,n}$ tends
to zero. Thus the hypothesis \eqref{eq:abstract-L2} of
Theorem~\ref{thm:abstract-extreme} holds with
$\varepsilon_{k,n}=4C_{k,n}/\lambda_{k,n}$.
Theorem~\ref{thm:abstract-extreme} now proves the proposition.
\end{proof}

\begin{rem}[Quadratic invariance versus invariance of means]
The vanishing in \eqref{eq:abstract-L2} is essential.
The weaker identity
\[
 \int_{K_n}f(hu)\,d\mu_n(u)
 =\int_{K_n}f(u)\,d\mu_n(u)
 \qquad(h\in K_k)
\]
is automatic from Haar invariance and provides no fixed-point
information by itself. The criterion instead controls the mean
square of the difference between an observable and its translate.
\end{rem}

\subsection{Application to the trace-class special unitary group}
\label{subsec:extreme-application}

\paragraph{Extension of the spectral estimates to Lipschitz functions.}
To apply Proposition~\ref{prop:abstract-spectral}, we make explicit
the extension of Lemmas~\ref{lem:relative-gap} and
\ref{lem:cross-energy-bound} beyond smooth functions.
For a real-valued $d_{1,n}$-Lipschitz function $f$ on $\SU(n)$,
define $\mathcal E_{k,n}(f)$ by \eqref{eq:cross-energy},
interpreting the derivatives $D_{X_{ij}}f$ and $D_{Y_{ij}}f$
in the weak sense. Since $n$ is fixed, the trace metric is locally
equivalent to a smooth Riemannian metric on $\SU(n)$;
thus $f\in W^{1,2}(\SU(n))$.
Convolution on the right with a smooth approximate identity gives
smooth functions $f_m$ converging to $f$ both uniformly and in
$W^{1,2}$. The metric $d_{1,n}$ is right-invariant, so
$\operatorname{Lip}_{d_{1,n}}(f_m)\le
\operatorname{Lip}_{d_{1,n}}(f)$.
The boundedness of $P_k$ on $L^2$ and the $W^{1,2}$ convergence
allow us to pass to the limit in both lemmas. Consequently,
for every such Lipschitz function $f$,
\begin{equation}\label{eq:extended-spectral-bounds}
\begin{aligned}
 (n-k)\|f-P_kf\|_{L^2(\nu_n)}^2
 &\le \mathcal E_{k,n}(f),\\
 \mathcal E_{k,n}(f)
 &\le 2k\operatorname{Lip}_{d_{1,n}}(f)^2.
\end{aligned}
\end{equation}

In particular, the proof of
Proposition~\ref{prop:abstract-spectral} gives
\begin{equation}\label{eq:asymptotic-L2}
\sup_{h\in\SU(k)}
\int_{\SU(n)}
|f(hu)-f(u)|^2\,d\nu_n(u)
\le
\frac{8k}{n-k}
\operatorname{Lip}_{d_{1,n}}(f)^2
\end{equation}
for every real-valued Lipschitz function
$f$ on $\SU(n)$.

\begin{proof}[Proof of Theorem~\ref{thm:SU1-extreme}]
Take $G=\SU_1(H)$, $K_n=\SU(n)$ and $d=d_1$.
The trace metric is bi-invariant, hence right-invariant, and
Lemma~\ref{lem:SU1-density} gives the required density.
For each $2\le k<n$, the projection $P_{k,n}$ in
Proposition~\ref{prop:abstract-spectral} is precisely $P_k$.
By \eqref{eq:extended-spectral-bounds}, the hypotheses
\eqref{eq:abstract-gap} and \eqref{eq:abstract-energy} hold
for all real-valued Lipschitz functions with
\[
  \lambda_{k,n}=n-k,\qquad C_{k,n}=2k.
\]
For each fixed $k\ge2$,
\[
  \frac{C_{k,n}}{\lambda_{k,n}}
  =\frac{2k}{n-k}\longrightarrow0
  \qquad(n\to\infty).
\]
Therefore Proposition~\ref{prop:abstract-spectral}, which was
deduced from Theorem~\ref{thm:abstract-extreme}, applies directly
and proves that $\SU_1(H)$ is extremely amenable.
\end{proof}

\begin{rem}
The proof does not show that $\SU(n)$ is a L\'evy family for
the trace metric; in fact this is ruled out by
Theorem~\ref{thm:compact-uniform}. It establishes the weaker but
fixed-point-sufficient property
\eqref{eq:asymptotic-L2} for every fixed $\SU(k)$.
Let $G_{\mathrm{ind}}$ denote the union
$\bigcup_n\SU(n)$ equipped with the direct-limit
topology in the category of topological groups, namely,
the finest group topology making all the inclusions
$\SU(n)\to G_{\mathrm{ind}}$ continuous.
The natural homomorphism
$G_{\mathrm{ind}}\to\SU_1(H)$ is continuous and has
dense image.

The finite-dimensional estimate
\eqref{eq:asymptotic-L2} remains valid independently
of the topology on the union. However, the abstract fixed-point
  criterion requires trace-norm uniformly continuous
orbit functions, which can be approximated by Lipschitz
functions with constants independent of $n$.
Continuity with respect to the direct-limit topology
does not provide this property. Thus our argument
does not settle the extreme amenability of
$G_{\mathrm{ind}}$.
\end{rem}

\section{The L\'evy property of Schatten-class unitary groups}
\label{sec:schatten-levy}

As observed by Pestov
\cite[Remark~4.1.12]{Pestov},
$\bigcup_{n\ge2}\SU(n)$ is dense in $\U_p(H)$
for $p>1$. We prove Proposition~\ref{prop:Up-levy}
using the concentration of Haar measures on $\SU(n)$.

\begin{proof}[Proof of Proposition~\ref{prop:Up-levy}]
Since the case $p=1$ is already known, as recalled in the
Introduction, it suffices to prove that $\U_p(H)$ is a
L\'evy group for every $p>1$.
Equip $\SU(n)$ with normalized Haar measure $\nu_n$ and the
Hilbert--Schmidt distance
$d_{2,n}(u,v):=\|u-v\|_{S_2}$.
We equip $\SU(n)$ with the bi-invariant Riemannian
metric determined by
\[
  g_e(X,Y):=-\Tr(XY),
  \qquad X,Y\in\mathfrak{su}(n).
\]
The Killing form of $\mathfrak{su}(n)$ is
\[
  B(X,Y)=2n\Tr(XY).
\]
Since the Ricci tensor of a bi-invariant Riemannian
metric on a compact Lie group is $-B/4$, we obtain
\[
  \mathrm{Ric}=\frac{n}{2} g.
\]
The Ricci-curvature concentration theorem of Gromov and Milman
\cite{GromovMilman}, applied to this metric, therefore yields
constants $c,C>0$, independent of $n$, such that
\begin{equation}\label{eq:GM-SU-concentration}
  \alpha_{(\SU(n),d_{2,n},\nu_n)}(\varepsilon)
  \le C\exp(-c n\varepsilon^2),
  \qquad \varepsilon>0.
\end{equation}
Here the Riemannian geodesic distance dominates $d_{2,n}$,
so the same bound holds for the latter. In particular,
$\bigl(\SU(n),d_{2,n},\nu_n\bigr)_{n\ge2}$ is a L\'evy
family. Later estimates for Haar unitary matrices provide another
route to such bounds \cite{MeckesMeckes}.

For $u,v\in\SU(n)$, define the restricted Schatten distance by
\[
  d_{p,n}(u,v):=\|u-v\|_{S_p}.
\]
The elementary comparison of singular-value norms then gives
\[
 d_{p,n}(u,v)\le
\begin{cases}
   n^{1/p-1/2}d_{2,n}(u,v),&1\le p<2,\\
   d_{2,n}(u,v),&2\le p<\infty.
\end{cases}
\]
Consequently, a comparison of metric neighborhoods in
\eqref{eq:GM-SU-concentration} gives
\begin{equation}\label{eq:p-concentration-background}
  \alpha_{(\SU(n),d_{p,n},\nu_n)}(\varepsilon)
  \le C\exp\!\left(-c\varepsilon^2
       n^{\min\{1,\,2-2/p\}}\right).
\end{equation}
The exponent of $n$ is strictly positive when $p>1$,
so the finite-dimensional special unitary groups form a L\'evy
family for each such $p$. By the density recalled above,
$\U_p(H)$ is a L\'evy group for all $p>1$.
\end{proof}

\section*{Acknowledgement}
ChatGPT was used in preparing this manuscript. The authors are responsible
for its mathematical content.

\end{document}